\input amstex
\documentstyle{amsppt}
\ifx\undefined\headlineheight\else\pageheight{47pc}\vcorrection{-1pc}\fi
\nologo

\font\scten= cmcsc10
\font\scten=cmcsc10

\loadeusm

\def\ga{\alpha}
\def\gb{\beta}
\def\gm{\gamma}
\def\gl{\lambda}

\def\gs{\sigma}
\def\gep{\epsilon}

\def\cv{\operatorname{conv}}

\def\conv{\operatorname{conv}}

\def\vvec{\vec v}
\def\wv{\vec w}
\def\xv{\vec x}
\def\vv{\vec v}
\def\yv{\vec y}
\def\zv{\vec z}

\def\rnn#1{\Bbb R^#1}

\def\lan{\langle}
\def\ran{\rangle}
\TagsOnRight
\document
\topmatter
\title
Locally Minimal Projections in  $B(H)$
\endtitle
\author
Charles A. Akemann and Joel Anderson
\endauthor
\thanks
Each author was partially supported by the National Science Foundation
during the period of
research that resulted in this paper.
\endthanks
\address
Department of Mathematics,
University of California,
Santa Barbara, CA 93106
\endaddress
\email
akemann\@math.ucsb.edu
\endemail
\address
Department of Mathematics,
Pennsylvania State University,
University Park, PA 16802
\endaddress
\email
anderson\@math.psu.edu
\endemail
\abstract
 Given an $n$-tuple $\{a_1, ..., a_n\}$ of self-adjoint
operators
 on an infinite dimensional Hilbert space $H$,
we say that a projection  $p$  in  $B(H)$ is {\bf locally minimal}
for  $\{a_1, ..., a_n\}$ if each
$pa_jp$ (for $j = 1, ..., n$) is a scalar multiple of $p$.  In Theorem 1.8 we
show that for any such $\{a_1, ..., a_n\}$ and any positive integer $k$
there
exists a projection $p$ of rank  $k$ that is locally minimal for $\{a_1, ...,
a_n\}$. If we further assume that
$\{a_1,\dots,a_n,1\}$ is a linearly independent set in the Calkin algebra,
then in Theorem 2.8 we prove that
$p$ can be chosen of infinite rank.
\endabstract
\endtopmatter
\document
\bigskip
\noindent {\scten \S0 Introduction } 

Given an $n$-tuple
$\{a_1, ..., a_n\}$ of bounded self-adjoint
operators on an infinite dimensional Hilbert space $H$,
we say that a projection  $p$  in  $B(H)$ (the algebra of all bounded
linear operators on  $H$)
is {\bf locally minimal} for  $\{a_1, ..., a_n\}$ if each
$pa_jp$ (for $j = 1, ..., n$) is a scalar multiple of $p$.  The name is
derived from the
fact that a projection $p$ in a C*-algebra  $A$  is a minimal projection if
$pAp$ consists
of multiples of  $p$.  In Theorem 1.8 we
show that
for any such $\{a_1, ..., a_n\}$ and any positive integer
$k$ there
exists a projection $p$ of rank  $k$ that is locally minimal for $\{a_1, ...,
a_n\}$.  If we further assume that
$\{a_1,\dots,a_n,1\}$ is a linearly independent set in the Calkin algebra
(i.e. the quotient of
 $B(H)$ by the ideal of compact operators on  $H$), then in
Theorem 2.8 we prove that
$p$ can be chosen of infinite rank.  In example 2.9 we show why our
theorems cannot be
improved, in general, although better results are often possible in
special cases, using our methods.

When only a finite set of operators is under study, we are not interested
in the algebra generated by these operators.  In this
situation local minimality appears to be the right concept to study.  For
instance, suppose
bounded self-adjoint operators
$\{a_1, ..., a_n\}$ represent physical
observables in a quantum system, states are represented by vectors in
$H$, and the expectation
value of an observable $a_j$ in a state
$\eta$ is given by
$\lan a_j\eta, \eta\ran$. If $p$ is locally minimal for $\{a_1, ...,
a_n\}$, then
as $\eta$
moves among the states under $p$ (i.e. such that $p\eta = \eta$), then the
expectation value of each of
the observables $\{a_1, ..., a_n\}$ does not change.  In a later paper we
plan to study
some physically interesting examples and see if a physical interpretation
can be imputed
to the phenomenon of local minimality.  Some of the more promising examples
will require
the methods of this paper to be extended to unbounded operators, where the
technical
difficulties are non-trivial.  For other physical applications it is
natural to assume
that the operators $\{a_1, ..., a_n\}$ lie in a type II or type III factor
and then ask
that the projection
$p$ lie in that same von Neumann algebra.  Our present methods do not
generalize to that
case.

When we began this study we wanted to carry over Samet's extension of
Lyapunov's
theorem \cite{3, p. 471} to non-commutative situations in the same way that we
generalized Lyapunov's theorem itself to non-commutative situations in our
Memoir
\cite{1}.  A key idea in Samet's proof was the notion of complete linear
independence.  A
set of functions on a measure space is {\bf completely linearly
independent} if, for any
set
$E$ of positive measure, the restriction of the set of functions to $E$ is
linearly
independent. Samet shows that complete linear independence is the generic
case.  As
we show in our Theorems 1.8 and 2.8, the natural non-commutative version of
complete
linear independence fails in  $B(H)$.  Therefore a generalization of Samet's
theorem will have to come from other methods.

\bigskip
\noindent\S 1 {\scten The Finite Rank Theorem}
\medskip
\proclaim{Lemma 1.1}  If
$$
a = \bmatrix
\ga_1 &0\\
0&\ga_2
\endbmatrix
\text{ and }
b =
\bmatrix
\gb_1&0\\
0&\gb_2
\endbmatrix
$$
are $2\times 2$ diagonal self--adjoint matrices such that
$$
\ga_1  \ge \gb_1\quad\text{and}\quad \gb_2 \ge \ga_2,
$$
then there is a vector $\eta$ in $\Bbb C^2$ such that
$$
\lan a\eta,\eta\ran = \lan b\eta,\eta\ran.
$$
\endproclaim
\demo{Proof} We may assume that $\ga_1 > \gb_1$ and $\gb_2 > \ga_2$ since
otherwise, the
assertion is trivial. Write
$$
t =\frac{\gb_2-\ga_2}{\ga_1-\gb_1 + \gb_2-\ga_2}.
$$
and set
$$
\eta_t =
\bmatrix
\sqrt{t}\\
\sqrt{1-t}
\endbmatrix.
$$
Note that since $\ga_1 > \gb_1$ and $\ga_2 < \gb_2$, we have $0 < t < 1$ and so
$\|\eta_t \| = 1$.  It is easy to check that, with this choice, we have
$$
\lan a\eta_t,\eta_t\ran =\lan b\eta_t,\eta_t\ran.
$$
\qed
\enddemo

\proclaim{Lemma 1.2} Suppose $k$ is a positive integer,  $a,b$,
and $c_1,\dots,c_m$
are diagonal self--adjoint matrices of dimension at least $2k$ and
 $\ga_1,\ga_2,\dots \ga_{2k} $ and $\gb_1,\gb_2\dots,\gb_{2k}$, are the
eigenvalues
lists for $a$ and $b$, respectively. If  $\ga_i \ge \gb_i$ for $i =
1,2,\dots, k$ and
$\ga_i \le \gb_i$ for $i = k+1, k+2, \dots, 2k$, then there is a projection $p$
of rank $k$ whose range is contained in the span of the first $2k$
eigenvectors of $a$
(or $b$) such that
\roster
\item $pap = pbp$.
\item $pap,pbp$ and $pc_1p,\dots,pc_mp$ can be simultaneously diagonalized.
\endroster
\endproclaim
\demo{Proof}  Write $\eta_1,\eta_2, \dots, \eta_k,\eta_{k+1},\dots,\eta_{2k}$
for the eigenvectors corresponding to the eigenvalues $\ga_1,$
$\ga_2,\dots,\ga_k,\ga_{k+1},\dots,\ga_{2k}$ and consider the
$2\times 2$  matrices $a_i$ and $b_i$ formed by
\noindent restricting $a$ and $b$ to the
span of
$\{\eta_i,\eta_{k+i}\}$, $i = 1,2,\dots,k$.  These matrices have the form
$$
a_i =
\bmatrix
\ga_i&0\\
0&\ga_{k+i}
\endbmatrix
\text{ and }
b_i = 
\bmatrix
\gb_i&0\\
0&\gb_{k+i}
\endbmatrix.
$$
Since $\ga_i \ge \gb_i$ and $\ga_{k+i} \le \gb_{k+i}$, we may apply  Lemma
1.1 to get
 a unit vector $\xi_i$ in the span of $\{\eta_i,\eta_{k+i}\}$
such
that
$(a\xi_i,\xi_i) = (b\xi_i,\xi_i)$. Let $p$ denote the projection onto the
span of
$\{\xi_1,\dots,\xi_k\}$.  Since distinct $\xi_i$'s are formed from orthogonal
eigenvectors we get that $\{\xi_1,\dots,\xi_k\}$ is orthonormal and  if $i\ne
j$, then
$(a\xi_i,\xi_j) = (b\xi_i,\xi_j)=(c_1\xi_i,\xi_j) =\cdots = (c_m\xi_i,\xi_j) =
0$.  Since
$(a\xi_i,\xi_i) = (b\xi_i,\xi_i)$ for each $i$, we  get $pap = pbp$
and $pap,pbp, pc_1p,\dots,pc_mp$ are simultaneously diagonalizable, as desired.
\qed
\enddemo

\proclaim{Lemma 1.3}  If $k$ is a positive integer and  $\ga_1,\dots,\ga_{3k}$
and $\gb_1\dots,\gb_{3k}$ are sequences of positive numbers,  then either
\roster
\item  $\{i:\ga_i = \gb_i\}$ has at least $k$ elements or,
\item  there is  an index set $\{i_1,\dots,
i_{2k}\}\subset \{1,2,\dots,3k\}$ and a positive real number $t$ such that
$$
\align
&\ga_{i_1} \ge t\gb_{i_1}, \dots, \ga_{i_k} \ge t\gb_{i_k} \text{ and }\\
&\ga_{i_{k+1}} \le t\gb_{i_{k+1}}, \dots, \ga_{i_{2k}} \le t\gb_{i_{2k}}.
\endalign
$$
\endroster
\endproclaim
\demo{Proof}If $\{i:\ga_i = \gb_i\}$ has at least $k$ members, we are done.  So
suppose that this set has fewer than $k$ elements.  In this case, relabeling
if necessary, we may assume that $\ga_i \ne \gb_i$ for $i = 1,2,\dots,2k$.  Now
write $\gm_i = \ga_i/\gb_i$ for $i=1,2,\dots, 2k$ and relabel so that $\gm_1
\ge \gm_2\ge \cdots \ge \gm_{2k}$. Observe that each $\gm_i > 0$ since the
$\ga$'s and
$\gb$'s are positive.   Now take
$t =
\gm_k$  With this we have
$$
\frac{\gm_1}{t} \ge \frac{\gm_2}{t}  \cdots \ge \frac{\gm_k}{t} = 1 \ge
\frac{\gm_{k+1}}{t} \ge \frac{\gm_{k+2}}{t} \ge \cdots \ge \frac{\gm_{2k}}{t}.
$$
In other words in this case we get
$$
\align
&\ga_1 \ge t\gb_1,\phantom{.} \ga_2 \ge t\gb_2,\phantom{.}\dots,\phantom{.}
\ga_k = t\gb_k\\
&\ga_{k+1} \le t\gb_{k+1},\phantom{.} \ga_{k+2} \le t\gb_{k+2 },\phantom{.}
\dots,\phantom{.}\ga_{2k} \le t\gb_{2k},
\endalign
$$
as desired.\qed
\enddemo

\proclaim{Lemma 1.4} If $a_1,a_2,\dots,a_n$ are self--adjoint operators
acting on
an infinite dimensional Hilbert space $H$, then there is a projection $p$ of
infinite rank such that $pa_1p,pa_2p,\dots, pa_np$ are each diagonal with
respect to some fixed orthonormal basis for the range of $p$.
\endproclaim

\demo{Proof} Fix a unit vector $\xi_1$ in $H$ and select a unit vector,
$\xi_2$ such that
$$
\xi_2 \in \{\xi_1,a_1\xi_1,\dots, a_n\xi_1\}^\perp
$$
Let us now continue by induction.  Suppose that orthonormal vectors
$\xi_1,\dots,\xi_k$
have been selected such that if $i < j \le k$, then
$$
\xi_j\in \{\xi_i,a_1\xi_i,\dots, a_n\xi_i:1 \le i < j\}^\perp.
$$
As $H$ is infinite dimensional, we may select a unit vector $\xi_{k+1}$
such that
$$
\xi_{n+1}\in \{\xi_i,a_1\xi_i,\dots, a_n\xi_i:1 \le i \le n\}^\perp.
$$

This produces an infinite orthonormal sequence $\xi_1,\xi_2,\dots$ such that
$$
\lan a_i\xi_j,\xi_k\ran = 0 \text{ if } j < k \text{ and } 1 \le i \le n.
$$
Since each $a_i$ is self-adjoint we also get $\lan a_i\xi_j,\xi_k\ran = 0
\text{ if } j >
k $.  Hence, if we let $p$ denote the projection onto the span of
$\xi_n$'s, then each
$pa_ip$ is diagonal.
\qed
\enddemo

\proclaim{Lemma 1.5} If $a_1,\dots, a_n$ are injective
operators acting on an infinite dimensional Hilbert space such that each
$a_i$ is either
positive or negative and $k$ is a positive integer, then  there is a
projection $p$ of
rank $k$ and nonzero real numbers
$t_2,\dots,t_n$ such that
$$
pa_1p = t_2pa_2p = \cdots = t_npa_np.
$$
\endproclaim
\demo{Proof} First assume that each $a_i$ is positive and injective.  By
Lemma 1.4
we may find a
projection
$p_1$ of infinite rank such that
$p_1a_ip_1$ is diagonal for each $i$.

Now write $N = 3^{n-1}k$.  If $N/3$ of the eigenvalues of the first $N$
eigenvalues of $p_1a_1p_1$ agree with the corresponding eigenvalues of
$p_1a_2p_1$, then
let
$p_2 <p_1$ denote the projection onto
the span of the associated eigenvectors.  Otherwise, by Lemma 1.3 we may
find a number
$t_2 > 0$ such that $N/3$ of the first $N$ eigenvalues of $p_1a_1p_1$ are
greater
than or equal to
the corresponding eigenvalues of
$t_2p_1a_2p_1$ and also find $N/3$ disjoint indices among the first $N$
eigenvalues where
this inequality is reversed.
Next, applying Lemma 1.2 we may find a projection $p_2 < p_1$ of rank $N/3$
such that
$$
p_2a_1p_2 = t_2p_2a_2p_2
$$
and such that $p_2a_1p_2,\dots,p_2a_np_2$ are all diagonal.

Now proceed by induction.  At the next stage we apply the argument of the
previous paragraph to get a projection $p_3 \le p_2$ of rank $N/3^2$   and
a positive
real number $t_3$ such that
$$
p_3a_1p_3 = t_3p_3a_3p_3
$$
and $p_3a_1p_3,\dots,p_3a_np_3$ are all diagonal.

 After $n - 1$ steps we get a projection $p = p_{n-1}$
of rank $N/3^{n-1} = k$ and positive numbers $t_2,\dots,t_n$ such that
$$
pap = t_2pa_2p = \cdots = t_npa_np.
$$

Thus the Lemma has been established  when all of the operators are
positive.  Now suppose
that some of the operators are negative.  Relabeling if necessary, we may
assume that for
some $1 \le k < n$ $a_1,\dots,a_k$ are positive and $a_{k+1},\dots, a_n$
are negative.
In this case, we may apply the first part of the proof to
$a_1,\dots,a_k,-a_{k+1},\dots,-a_n$ to get a projection $p$ of rank $k$ and
scalars
$t_2,\dots,t_n$ such that
$$
pa_1p = t_2pa_2p = \cdots = t_kpa_kp = -t_{k+1}pa_{k+1} = \cdots = -t_npa_np.
$$
The proof may now be completed by replacing $t_{k+1},\dots, t_n$ with
$-t_{k+1},\dots,
-t_n$.
\qed
\enddemo

\proclaim{Lemma 1.6} If $a$ is a self-adjoint $2k\times 2k$ matrix, then
there are
orthonormal vectors $\xi_1,\dots,\xi_k$ in $\Bbb C^{2k}$such that if $p$
denotes the projection
onto the span of the $\xi_i$'s, then $pap = \gl p$ for some $\gl \in \Bbb R$.
\endproclaim
\demo{Proof} We may assume that $a$ is diagonal with diagonal entries
$\ga_1 \le \dots \le
\ga_{2k}$. Let $\eta_1,\dots,\eta_{2k}$ denote the standard basis vectors
and fix a real
number $\gl$ such that
$$
\ga_k \le \gl \le \ga_{k+1}.
$$
Next if $i \le k$, then write
$$
t_i = \cases  0 &\text{ if } \ga_{i} = \gl = \ga_{k+i}\\
\dsize\frac{\gl - \ga_i}{\ga_{k+i}-\ga_i} &\text{ otherwise} .
\endcases
$$
If we  set
$$
\xi_i = \sqrt{t_i}\eta_i +
\sqrt{1-t_i}\eta_{k+i},\quad 1 \le i \le k,
$$
then it is straightforward to check that the $\xi_i$'s are orthonormal and
$\lan
a\xi_i,\xi_i\ran =
\gl$ for each index $i$.  Further, since $a$ is diagonal, we get that
$$
\lan a\xi_i,\xi_j\ran = 0 \text{ if } i \ne j.
$$

Hence, if $p$ denotes the projection onto the span of the $\xi_i$'s, then
$pap = \gl p$.
\qed
\enddemo

\proclaim{Lemma 1.7} If $n$ is a positive integer and for each $1 \le i \le
n$,
$\{a_{ik}\}_{k=1}^\infty$ is an infinite real sequence, then there is
an infinite index set $\{j_1,j_2,\dots\}$ such that for each $i$ the
subsequence $\{\ga_{i,j_1},\ga_{i,j_2},\dots\}$ is either strictly
positive, strictly
negative or 0.
\endproclaim
\demo{Proof} Write
$$
\align
\gs_{1,+}&= \{\ga_{1j}: \ga_{ij} > 0\},\\
\gs_{1,-} &= \{\ga_{1j}:\ga_{1j} < 0\},\\
\gs_{1,0}& = \{\ga_{1j}: \ga_{1j} = 0\}
\endalign
$$
and let $\gs_1$ denote the first of these sets that is infinite.  Now
continue by
induction.  Suppose that for some index $1 \le i < n$, infinite subsets
$\gs_1\supset
\cdots \supset \gs_i$ have been selected such that  if $1 \le j \le i$,
then $\{\ga_{jk}:
k\in \gs_j\}$ is either strictly positive, strictly negative or 0.
Arguing as above, we
may find an infinite subset $\gs_{i+1}$ of $\gs_i$ such  that $\{\ga_{i+1,k}:
k\in \gs_{i+1}\}$ is either strictly positive, strictly negative or 0.  The
proof is now
completed by taking
$$
\{j_1,j_2\dots\} = \gs_n.
$$
\qed
\enddemo

\proclaim{Theorem 1.8 (The Finite Rank Theorem)} If $a_1,\dots, a_n$ are
bounded,  self--adjoint operators acting on an infinite dimensional Hilbert
space $H$ and
$k$ is a positive integer, then there are real numbers $t_2,\dots,t_n$ and
a projection
$p$ of rank $k$ such that
$$
pa_ip = t_ip, \qquad 1 \le i \le n.
$$
\endproclaim
\demo{Proof}  Applying Lemma 1.4, we may  find a projection $q_1$ of
infinite rank such
that each $q_1a_iq_1$  is diagonal with respect to some fixed orthonormal
basis for
$q_1H$.  Now apply Lemma 1.7 to the eigenvalue sequences of the
$q_1a_iq_1$'s to get a
diagonal projection $q_2 \le q_1$ of infinite rank and such that for each
index $i$,
$q_2a_iq_2$ is either strictly positive, strictly negative or 0 on $q_2H$.

Let $t_i = 0$ for each index $i$ such that  $q_2a_iq_2 = 0$.  If $q_2a_iq_2
= 0$ for each
$i$, then we may take $p = q_2$ and the proof is complete.  Otherwise,
restricting to the
nonzero  $q_2a_iq_2$'s, we may assume that  $q_2a_iq_2$ is injective for
each $i$.
In this case we may apply Lemma 1.5 and  get a projection $p_1
\le q_2$ of rank at least $2k$ and $t_2,\dots, t_{n} \ne  0$ such that
$p_1a_1p_1 =
t_2p_1a_2p_1 = \cdots = t_{n}p_1a_np_1$.  The proof is completed by
applying Lemma 1.6
to $p_1a_1p_1$ to get a projection $p$ of rank $k$ with the desired properties.
\qed
\enddemo

\bigskip
\noindent \S2 {\scten The Infinite Rank Theorem}

\definition{Definition 2.1} If $x_1,x_2,\dots,x_n$ are bounded linear operators
on an infinite dimensional separable Hilbert space $H$, their {\bf joint
essential numerical range} is by definition
$$
W_e(x_1,x_2,\dots,x_n) = \left\{\left(f(x_1),f(x_2),\dots,f(x_n)\right) :f
\text{
is a singular state on } B(H)\right\}.
$$
\enddefinition

\definition{Remark}
Recall that a singular state of  B(H) is a positive linear functional of norm 1 whose
kernel contains the
compact operators on  $H$.  By the Hahn-Banach Theorem, an operator in
$B(H)$ that is in the
kernel of all singular states must be a compact operator.  This fact will
be used in Lemma 2.7.
Note that
$W_e(x_1,x_2,\dots,x_n)$  is a convex subset of
$\Bbb C^n$ and, if each $x_i$ is self-adjoint, then
$W_e(x_1,x_2,\dots,x_n)\subset \Bbb R^n$.
\enddefinition

\definition{Notation}
 If $x_1,\dots,x_n$ are operators on a  Hilbert space $H$ and $\eta
\in H$, we write
$$
\vvec(x_1,\dots,x_n,\eta) = (\lan x_1\eta,\eta\ran,\dots, \lan
x_n\eta,\eta\ran)
$$
for the vector in $\Bbb C^n$ (or $\Bbb R^n$ if the $x_i$'s are
self-adjoint) which these
elements determine.
\enddefinition

\proclaim{Lemma 2.2} With notation as above, if $\vvec =
(\gl_1,\gl_2,\dots,\gl_n) \in
W_e(x_1,x_2,\dots,x_n)$,
$F$ is a finite orthonormal set of vectors in $H$ and $\epsilon > 0$, then
there
is a unit vector $\eta \in H$ such that
\roster
\item $\eta$ is orthogonal to each vector in $F$.
\item $\|\vvec
-\vvec(x_1,\dots,x_n,\eta)\| < \epsilon$.
\endroster
\endproclaim
\demo{Proof}   By definition of the
joint essential numerical range, there is a singular state  $f$  such that
$f(x_i) = \gl_i$  for  $i = 1, ..., n$. Let  $p$  denote the projection onto  $F^\perp$.  Since  $F$  is finite,
and $f$ is a singular state,  $f(pap) = f(a)$  for each  $a$  in  $B(H)$, so we may
regard $f$ as a singular state on $B(H)$. By a theorem of  Glimm \cite{2, Theorem 2, p.
216} there is a sequence  $\{\xi_k\}$  of
unit vectors in $pH$  such that
$$
\lim_{k\to\infty} \lan x_i\xi_k,\xi_k\ran = f(x_i) = \gl_i  \text{ for } i
= 1, ..., n.
$$
 Now select an integer $k$  such that
$$
|\lan x_i\xi_k, \xi_k\ran - f(x_i)| < \frac{\gep}{\sqrt{n}}
$$  for each  $i = 1, ..., n$,
and let  $\eta = \xi_k$.  Since  $\eta$ lies in  $pH$ , part $(1)$ of the
Lemma is satisfied.
Further, we have
$$
\align
\|\vv - \vv(x_1, ..., x_n,\eta)|| &= \|(\lan x_1\xi_k, \xi_k\ran, \dots, \lan
x_n\xi_k, \xi_k\ran) - (\gl_1,\dots,\gl_n)\| \\ & =  \|(\lan x_1\xi_k,
\xi_k\ran - f(x_1),
\dots, \lan x_n\xi_k, \xi_k\ran  - f(x_n))\|<
\gep
\endalign
$$
and so part $(2)$ of the Lemma is true.
\qed
\enddemo

\proclaim{Lemma 2.3} If $a_1,\dots,a_n$ are self-adjoint operators in $B(H)$,
$\vv_1,\dots,\vv_{n+1}$ are vectors in $W_e(a_1,\dots,a_n)$, $F_0$ is a finite
orthonormal set and $\epsilon > 0$, then there are orthonormal vectors
$\eta_{1},\dots,\eta_{n+1}$  and finite orthonormal sets $F_1,\dots F_{n+1}$ in
$H$  such that
\roster
\item $F_0 \subset F_1 \subset \cdots \subset F_{n+1}$,
\item $\eta_{i} \in F_{i} \cap F_{i-1}^\perp, \quad 1 \le i\le n+1$,
\item  $\{a_1\eta_{i},a_2\eta_{i},\dots,a_n\eta_{i} \}\subset
\text{span}(F_{i}),\quad 1 \le i \le n+1$,
\item $\|\vvec_i -\vvec(a_1,\dots,a_n,\eta_{i})\|< \gep$ for $i =
1,\dots,n+1$, and
\item if $j \ne k$, then $\lan a_i\eta_j,\eta_k\ran = 0,\quad 1 \le i \le
n,\quad 1 \le
j,k\le n+1$.
\endroster
\endproclaim
\demo{Proof} Applying Lemma  2.2  we may select $\eta_{1}$ in $H$ such that
$$
\| \vvec_1 - \vvec(a_1,\dots,a_n,\eta_{1})\|_2 < \gep
$$
and $\eta_1 \perp F_0$. Now select a finite orthonormal set $F_1$ such that
$F_0\cup\{\eta_1\} \subset
F_1$ and
$$
\{a_1\eta_{1},a_2\eta_{1},\dots,a_n\eta_{1} \}\subset
\text{span}(F_{1}).
$$

Continuing by induction,  suppose that vectors $\eta_1,\dots,\eta_i$ and
finite
orthonormal sets $F_1,\dots F_i$ satisfying conditions $(1),(2),(3)$ and
$(4)$  have been
selected for some $1 \le i < n+1$.  Applying Lemma 2.2 once more to $F_i$
we get a unit
vector $\eta_{i+1}$ such that $\eta_{i+1} \perp F_{i}$ and $\|\vvec_{i+1}
-\vvec(a_1,\dots,a_n,\eta_{i+1})\|< \gep$.  Now select a finite orthonormal
set $F_{i+1}$
containing $\eta_i$ and $F_i$ and such that
$$
\{a_1\eta_{i+1},a_2\eta_{i+1},\dots,a_n\eta_{i+1} \}\subset
\text{span}(F_{i+1}).
$$
This completes the induction argument.

Finally, observe that by construction if $j < k \le
n+1$, then
$\lan a_i\eta_j,\eta_k\ran = 0,\quad 1 \le i \le n$ and since each $a_i$ is
self-adjoint
these equalities also hold for $j > k$.  Hence, condition $(5)$ is true.
\qed
\enddemo

It is useful to record some facts about $n$--simplices before proceeding.
Recall that an {\bf$\bold n$--simplex }in $\Bbb R^n$ is a convex set of the
form
$$
S = \cv(\vv_1,\dots,\vv_{n+1})
$$
 such that its interior is not empty.  The vectors $\vv_1,\dots,\vv_{n+1}$
are called the
{\bf vertices} of $S$.  The following fact is well known and  easy to prove.

\proclaim{Lemma 2.4}  Every vector in an $n$--simplex has a unique
representation as a
convex combination of its  vertices.
\endproclaim

The coefficients in this
representation are called  the {\bf barycentric coordinates} of the vector.

Lemma 2.6 below is the key to the proof of the main Theorem in this section
(Theorem 2.8).
We are grateful to Nik Weaver for showing us an argument that considerably
shortens our
original proof.   We shall employ the following notation.  If $S$ is a
subset of $\Bbb
R^n$ and $\gep > 0$, then  the {\bf $\gep$-interior}  of $S$ is by
definition the set of
all vectors $\vv \in S$ such that $S$ contains the closed $\gep$-ball
centered at $\vv$.
The proof of the next result is routine.

\proclaim{Lemma 2.5} If $C$ is a convex subset of $\Bbb R^n$ and $\vv$ is a
vector in
the interior of $C$, then there is $\gep > 0$ and vectors  $\vv_1,\dots,
\vv_{n+1}$in $C$
such that
\roster
\item $S = \cv{(\vv_1,\dots,\vv_{n+1})} \subset C$ is an $n$-simplex and
\item $\vv$ is in the $\gep$-interior of $S$.

\endroster
\endproclaim

\proclaim{Lemma 2.6}    If
$$
S = \conv{(\vv_1,\dots,\vv_{n+1})}\text{ and } S^\prime =
\conv{(\vv_1^\prime,\dots
\vv_{n+1}^\prime)}
$$
are $n$-simplices in $\Bbb R^n$,  $\gep  > 0$, and
$$
\|\vv_i - \vv_i^\prime\| < \gep \text{ for }  i = 1,\dots, n+1,
$$
then $S^\prime$ contains the $\gep$-interior of $S$.   In particular, if
$\vv$ is in the
$\gep$-interior of $S$, then $\vv$ is in each  $S^\prime$ that  satisfies
the hypotheses.
\endproclaim
\demo{Proof}  Suppose that there is a vector $\vv$ in the $\gep$-interior
of $S$ such
that $\vv$ is not in $S^\prime$.  Since $S^\prime$ is convex, there is a
hyperplane $P$
that strictly separates $\vv$ and $S^\prime$.   Let $\xv$ denote the unique
vector in $P$
that is closest to $\vv$ and write $\yv = \vv -\xv$.  We have that $\yv$ is
orthogonal to
$P$ and ``points away'' from $S^\prime$.  Now set
$$
\wv = \vv + \frac{\gep}{\|\yv\|}\yv.
$$
Since $\|\wv-\vv\| = \gep$ and $S$ contains the closed $\gep$-ball about
$\vv$, $\wv \in
S$.  By construction, $\xv$ is the unique vector in $P$ that is closest to
$\wv$.  Thus,
if
$\zv\in P$, then  we have
$$
\|\wv-\zv\| \ge \|\wv - \xv\| =  \|\wv - \vv + \yv\| = \left( 1 +
\frac{\gep}{\|\yv\|}\right)\|\yv\| =
\|\yv\| + \gep > \gep
$$
and so $\text{dist}(\wv, S^\prime) > \gep$.

On the other hand, as $\wv \in S$, by Lemma 2.4 we have
$$
\wv = \sum_{i=1}^{n+1} t_i\vv_i,
$$
where the $t_i$'s are the barycentric coordinates of $\wv$.
If we now write
$$
\wv^\prime = \sum_{i=1}^mt_i\vv_i^\prime,
$$
then we have $\wv^\prime \in S^\prime$ and
$$
\|\wv -\wv^\prime\|   \le \sum_{i=1}^m t_i \|\vv_i - \vv_i^\prime\|   < \gep,
$$
which is a contradiction.
\qed
\enddemo

\proclaim{Lemma 2.7} If $a_1,\dots,a_n $ are bounded self--adjoint linear
operators acting on a separable Hilbert space $H$, then
$a_1,\dots,a_n, 1$ are linearly independent in the Calkin algebra
if and only if $W_e(a_1,\dots,a_n)$ has nonempty interior in
$\Bbb R^n$.
\endproclaim
\demo{Proof} Suppose that $W_e(a_1,\dots,a_n)$ has empty interior.  This means
that the dimension of $W_e(a_1,\dots,a_n)$ is less than $n$ and so
$W_e(a_1,\dots,a_n)
\subset \vv + W$, where $\vv\in \rnn{n}$ and $W$ is a subspace of dimension
less than
$n$.  Let $\xv = (t_1,\dots,t_n)$ denote a unit vector that is orthogonal
to $W$.  If $f$ is a singular state on $B(H)$ so that
$$
(f(a_1),\dots,f(a_n) = \vv + \wv
$$
for some $\wv \in W$, and we write $t_{n+1} =
-\lan\vv,\xv\ran$, then we have
$$
\align
f(t_1a_1 +\cdots+t_na_n + t_{n+1}I) &= \lan(f(a_1),\dots,f(a_n)),\xv\ran
-\lan\vv,\xv\ran \\
& = \lan \vv + \wv,\xv\ran -\lan \vv,\xv\ran = 0
\endalign
$$
and therefore, as noted in  the Remark following Definition
2.1,  $t_1a_1 +\cdots + t_na_n + t_{n+1}I$ is
compact .

Conversely, if there are real scalars $t_1,\dots, t_n$ and $t$ such that
$$
t_1a_1 +\cdots + t_na_n + t1
$$
is compact, then  for every singular state $f$ we have
$$
\sum_{i=1}^n t_if(a_i) = -t
$$
and so $W_e(a_1,\dots,a_n)$ lies in an hyperplane.  Thus this set has empty
interior.
\qed
\enddemo

\proclaim{Theorem 2.8 (The Infinite Rank Theorem)} If $a_1,\dots, a_n$ are
bounded self--adjoint linear operators on a separable infinite dimensional
Hilbert space
$H$ such that $a_1,\dots,a_n$ and $1$ are linearly independent in the
Calkin algebra and
$$
\vv = (\gl_1,\dots,\gl_n)
$$
is a vector in the interior of $W(a_1,\dots,a_n)$, then there is a
projection $p$  of
infinite rank  on $H$ such that
$$
pa_ip = \gl_ip
$$
for each index $i$.
\endproclaim
\demo{Proof}  As $\vv$ is in the interior of $W(a_1,\dots,a_n)$,
 we may apply Lemma 2.5 to find $\gep > 0$ and vectors
$\vv_1,\dots,\vv_{n+1}$ in
$W(a_1,\dots,a_n)$ such that if $S = \cv{(\vv_1,\dots,\vv_{n+1})}$, then
$S$ is an
$n$-simplex, $S \subset W(a_1,\dots,a_n)$  and $\vv$ is in the
$\gep$-interior of $S$.

We shall now apply Lemma 2.3 in an inductive construction.  In each application
we shall use the operators $a_1,\dots,a_n$,  the vectors
$\vv_1,\dots,\vv_{n+1}$ and the number $\gep$  found in the previous
paragraph.  The finite orthonormal set used in each application will be
constructed by the inductive process.

Let  us begin by setting $F_0 = \emptyset$ and applying  Lemma 2.3 to get
orthonormal
vectors
$\eta_{11},\dots,\eta_{1,n+1}$ and finite orthonormal sets
$F_{11},\dots,F_{1,n+1}$
satisfying conditions $(1),(2),(3),(4)$ and $(5)$  of this Lemma.

Now let us proceed by induction. Suppose that for some integer $k \ge 1$,
sequences
$$
\eta_{11},\dots,\eta_{1,n+1},\dots, \eta_{k1},\dots,\eta_{k,n+1}
$$
and finite orthonormal sets

$$
F_{11} \subset \cdots \subset F_{1,n+1} \subset \cdots \subset
F_{k1}\subset \cdots \subset F_{k,n+1}
$$
have been selected such that if $1 \le j < k$,
then conditions $(1),(2),(3),(4)$ and $(5)$
of Lemma 2.3 hold for $F_{j,n+1}$, $F_{j+1,1},\dots, F_{j+1,n+1}$ and
$\eta_{j+1,1},\dots,
\eta_{j+1,n+1}$.

To complete the inductive argument, apply Lemma 2.3 once more to
$F_{k,n+1}$ to get
$F_{k+1,1},\dots, F_{k+1,n+1}$ and $\eta_{k+1,1},\dots,\eta_{k+1,n+1}$
satisfying
conditions $(1),(2),(3),$ $(4)$ and $(5)$ of Lemma 2.3.

This argument produces  an infinite orthonormal sequence
$$
\eta_{11},\dots,\eta_{1,n+1},\eta_{21},\dots,\eta_{2,n+1},\dots,\eta_{k1},\dots,
\eta_{k,n+1},\dots
$$
such that
$$
\|\vvec_i - \vvec(a_1,\dots,a_n,\eta_{ki})\| < \gep \quad i = 1,2,\dots,
,n+1,\,\, k
= 1,2,
\dots
$$
and
$$
\lan a_i\eta_{jk},\eta_{lm}\ran = 0  \quad 1 \le i \le n, \text{ if } j \ne
l \text{ or }
k \ne m.\tag$*$
$$

Next, for each index $k$ write
$$
S_k =\cv( \vvec(a_1,\dots,a_n,\eta_{k1}),\vvec(a_1,\dots,a_n,\eta_{k2}),\dots
,\vvec(a_1,\dots,a_n,\eta_{k,n+1})).
$$
Since $\|\vv_i-\vv(a_1,\dots,a_n,\eta_{ki})\| < \gep$ for each $i$ and
$\vv$ is in the
$\gep$-interior of $S$,  we have
$\vv \in S_k$ for each $k$ by Lemma 2.6. Let
$t_1,\dots,t_{n+1}$ denote the barycentric coordinates of for $\vv$ as an
element of $S$
and write $s_{k1},\dots,s_{k,n+1}$ for the barycentric coordinates of $\vv$
as an element
of
$S_k$  so that we have
$$
\vv = \sum_{i=1}^{n+1}t_i\vv_i =
\sum_{i=1}^{n+1}s_{ki}\vv(a_1,\dots,a_n,\eta_{ki})
$$
Now write
$$
\phi_k = \sum_{i=1}^{n+1}\sqrt{s_{ki}}\eta_{ki},\qquad k = 1,2,\dots.
$$
We have that $\{\phi_k\}$ is an infinite orthonormal set by construction.  Next
 let $p$ denote the projection onto the span of
$\phi_1,\phi_2,\dots$.  Observe that by $(*)$ above, if $j \ne k$, then
$$
(a_r\phi_k,\phi_j)  = (a_r\phi_j,\phi_k) = 0 \quad r = 1,\dots, n.
$$
Also, we have
$$
(a_r\phi_k,\phi_k) = \sum_{i=1}^{n+1} s_{ki}(a_r\eta_{ki},\eta_{ki})
$$
because by $(*)$ above the cross terms are $0$ and this means that
$$
\align
\vvec(a_1,\dots,a_n,\phi_k)&= ((a_1\phi_k,\phi_k),\dots,(a_n\phi_k,\phi_k)) \\
&=\bigg(\sum_{i=1}^{n+1}s_{ki}(a_1\eta_{ki},\eta_{ki}),\dots,
\sum_{i=1}^{n+1}s_{ki}(a_n\eta_{ki},\eta_{ki})\bigg)\\
&=\sum_{i=1}^{n+1}s_{ki}\vvec(a_1,\dots,a_n,\eta_{ki})\\
& =
\sum_{i=1}^{n+1}t_i\vvec_i = \vv = (\gl_1,\dots,\gl_n).
\endalign
$$
In other words, for each $k$,
$$
(a_r\phi_k,\phi_k) = \gl_r \quad r = 1,\dots,n+1.
$$
Hence,
$$
pa_ip = \gl_ip \quad\text{  for } i = 1,\dots,n.
$$
\qed
\enddemo
\noindent {\bf Examples.}
\roster
\item  Let $a$ denote any self-adjoint operator on $H$, suppose $b$
is a compact positive operator with dense range and write
$$
a_1 = a \text{ and } a_2 = a+b.
$$
Suppose that $p$ is a projection of infinite rank such that $pa_1p = tp$.
In this case we
have
$$
pa_2p = tp + pbp
$$
and so if $pa_2p$ is also a multiple of $p$, then $pbp$ is a multiple of
$p$.  Since $pbp$
is compact and $p$ has infinite rank, we get that $pbp = 0$.  But since $b$
is positive
and has dense range, this is impossible.  Thus no such projection exists.
\item  On the other hand, suppose $\{\eta_n\}$ denotes an orthonormal basis
for $H$ and
we define the compact operator $c$ by the formulas
$$
c\eta_{2n-1} = \frac{1}{n} \text{ and } c\eta_{2n} = -\frac{1}{n},\quad n =
1,2,\dots .
$$
If we then write
$$
\xi_n = \frac{1}{\sqrt{2}}(\eta_{2n-1}+\eta_{2n})
$$
and let $p$ denote the projection onto the span of the $\xi_n$'s, then we
have $pcp = 0$
and so the conclusion of Theorem 2.8 holds for $p$ and $p+c$.
\endroster

\bigskip
\centerline {\smc References} \Refs\nofrills{}
\bigskip

\noindent [1] Charles A.\ Akemann and Joel\ Anderson, Lyapunov Theorems for
Operator Algebras, {\it Mem.\ Amer.\ Math.\ Soc.\ \bf 458} (1991).
\medskip
\noindent [2] James\ Glimm, A Stone-Weierstrass Theorem for Operator Algebras,
{\it Annals of Math.\ \bf 72(2)} (1960), 216--244.
\medskip
\noindent [3] Dov\ Samet, Continuous Selections for Vector Measures,
{\it Mathematics\ of Operations Research.\ \bf 12} (1987), 536--543.

\smallskip
\endRefs

\enddocument